\documentclass[11pt]{amsart}
\usepackage{amsmath,amssymb,amsthm}
\usepackage{fullpage}
\usepackage{xcolor}
\newtheorem{theorem}{Theorem}[section]
\newtheorem{lemma}[theorem]{Lemma}
\newtheorem{proposition}[theorem]{Proposition}
\newtheorem{corollary}[theorem]{Corollary}
\theoremstyle{remark}
\newtheorem{remark}[theorem]{Remark}
\theoremstyle{definition}
\newtheorem{definition}[theorem]{Definition}
\newtheorem{example}[theorem]{Example}

\newcommand{\Hom}{\operatorname{Hom}}
\newcommand{\GL}{\operatorname{GL}}

\newcommand{\D}{\Delta}
\newcommand{\la}{\lambda}
\newcommand{\K}{\mathbb K}
\numberwithin{equation}{section}

\title[Homomorphism spaces for $\mathrm{GL}_3$ in characteristic $2$]{Homomorphisms between Weyl Modules for $\mathrm{GL}_3$ in Characteristic $2$}
\author[M. Maliakas]{Mihalis Maliakas}
\address{Department of Mathematics\\National and Kapodistrian University of Athens\\Athens, Greece}
\email{mmaliak@math.uoa.gr}

\author[D.-D. Stergiopoulou]{Dimitra-Dionysia Stergiopoulou}
\address{Department of Mathematics\\University of Thessaly\\Lamia, Greece}
\email{dstergiop@math.uoa.gr}
\email{dstergiopoulou@uth.gr}

\date{November 14, 2026}

\begin{document}

	\begin{abstract}
We classify the homomorphism spaces between Weyl modules for the general linear group $\GL_3$ over an algebraically closed field of characteristic $2$. We prove that these spaces have dimension at most one and give an explicit recursive procedure that determines when they are nonzero and constructs generators in terms of semistandard tableaux. As a consequence, we determine the dimensions of certain homomorphism spaces between Specht modules in characteristic $2$.
	\end{abstract}
	\maketitle

	\section{Introduction}\label{sec1}
	Let $G$ be a connected reductive algebraic group over an algebraically closed
	field $\K$ of positive characteristic $p$. The Weyl modules
	$\Delta(\la)$ are finite-dimensional rational $G$-modules indexed
	by dominant weights. Their structure is a central problem in the rational representation
	theory of $G$ and even the homomorphism spaces
	$\Hom_G(\Delta(\la),\Delta(\mu))$ are far from understood.

	For the general linear group $\GL_n$, nonzero homomorphisms $\Delta(\la) \to \Delta(\mu)$
	under suitable conditions on $\la, \mu$ and  $p$, where $\mu$ is obtained by raising cells from a row of $\la$, were obtained by Carter and Lusztig \cite{CL} and by Carter and Payne \cite{CPay}. Extensions of the Carter-Payne theorem were obtained by Fayers and Martin \cite{FM}. In the special case where a single cell is raised, it is known that the corresponding $\Hom$-space has dimension at most one; see Kulkarni \cite{Ku} and Ellers and Murray \cite{EM}. On the other hand, there are examples of $\Hom$-spaces that have dimension greater than one and even arbitrarily large dimensions, see Dodge \cite{Dodge} and Lyle \cite{Lyle}. Stability properties of their dimensions were studied in \cite{EMS}.

For $\mathrm{GL}_3$, the situation is better understood. When $\mu$ is obtained from $\la$ by moving cells from the third row to the first, the dimension of $\Hom_G(\Delta(\la),\Delta(\mu))$ follows from the computation of the corresponding integral $\operatorname{Ext}^1$ group \cite{FC,BF}. Cox and Parker \cite{CP} classified all homomorphisms
between Weyl modules for $\mathrm{GL}_3$ when $p>2$ and proved
that the  homomorphism spaces have dimension at most one. In this paper, we treat the remaining case $p=2$, thus completing the classification for $\GL_3$.

Our main result is a recursive classification of the homomorphism spaces between Weyl modules for $\GL_3$ in characteristic $2$ (Theorem \ref{thm2}). For any pair of dominant weights $\la, \mu$, the procedure determines whether $\Hom_G\bigl(\Delta(\la),\Delta(\mu)\bigr)$ is nonzero and, when it is, constructs an explicit generator in terms of semistandard tableaux. In particular, we have that
$\dim_\K\Hom_G\bigl(\Delta(\la),\Delta(\mu)\bigr)\le1$ for all $\la, \mu$ (Corollary \ref{dimhom}).

In Section \ref{sec2}, we discuss preliminary reductions and use results from \cite{MS} to identify the homomorphism spaces with the common kernels of two explicit matrices. In Section \ref{sec3}, we use Lucas' theorem to relate these kernels to those associated with smaller pairs of partitions. We also give explicit formulas for recovering the coefficients of a homomorphism from those of the smaller pair. In Section \ref{sec4}, we prove Theorem \ref{thm2}. In Section \ref{secspecht}, we combine our results with McDowell's theorem on the inverse Schur functor \cite{McD} to determine the dimensions of certain homomorphism spaces between Specht modules in characteristic $2$. This includes all pairs of partitions of $r\ge3$ for which the source has at most two parts and the target has at most three parts (Corollary \ref{specht}). Section \ref{lemproofs} contains the proofs of two lemmas used in Section \ref{sec3}.

	\section{Preliminaries}\label{sec2}

Throughout this paper, $\K$ denotes an algebraically closed field of characteristic $2$ and $G$ denotes the general linear group $\GL_3(\K)$.	We identify the weights of $G$ with triples $\alpha = (\alpha_1, \alpha_2, \alpha_3)$ of integers and the dominant weights of $G$ with triples $\la =(\la_1, \la_2, \la_3)$ of integers such that $\la_1 \ge \la_2 \ge \la_3$. For a dominant weight $\la$ we denote the corresponding Weyl module of $G$ by $\Delta(\la)$, see \cite[II.2.13]{Ja} where the notation $V(\la)$ is used. In the case where $\la$ is a partition, i.e. when $\la_1 \ge \la_2 \ge \la_3 \ge 0$, we use properties of $\D(\la)$ described in \cite[Section 2]{W}. In Weyman's book, $\D(\la)$ is denoted by $K_\la E$, where $E$ is the natural $\GL_3(\K)$-module.

\subsection{Reductions}\label{secred}If	$\alpha=(\alpha_1, \alpha_2, \alpha_3)$ and $\beta=(\beta_1, \beta_2, \beta_3)$ are weights,
	 we write $ \alpha \unlhd \beta$ if $\alpha_1 \le \beta_1$, $\alpha_1 +\alpha_2 \le \beta_1 +\beta_2$ and $\alpha_1 +\alpha_2 +\alpha_3 = \beta_1 +\beta_2 +\beta_3$. Since every weight $\alpha$ of $\Delta(\mu)$ satisfies $\alpha \unlhd \mu$, see \cite[II.2.2, Proposition (b)]{Ja} for the corresponding statement involving the duals of Weyl modules, we have that $\Hom_G(\Delta(\la),\D(\mu))=0$ unless $ \la \unlhd \mu$. Hence we may assume $\mu_3 \le \la_3$.

For an integer $m$ let $\det^{m}$ denote the $m$th power of the determinant representation. If $\la =(\la_1, \la_2, \la_3)$ is a dominant weight, then $\D(\la) \otimes \det^{m} \cong \D(\la_1+m, \la_2+m, \la_3+m)$. For any dominant weights $\la, \mu$ we have $\Hom_G(\D(\la), \D(\mu)) \cong \Hom_G(\D(\la) \otimes {\det}^{m}, \D(\mu) \otimes {\det}^{m})$. Taking $m=-\mu_3$, we may assume $\la$ and $\mu$ are partitions and $\mu$ has at most two parts.

	\begin{lemma}\label{lemtwist}
		Let $\la=(\la_1,\la_2,\la_3)$ and
		$\mu=(\mu_1,\mu_2)$ be partitions of $r$ with
		$\mu_2>\la_1$. Then
		$
		\Hom_G\bigl(\Delta(\la),
		\Delta(\mu)\bigr)
	$ is isomorphic to
	$\Hom_G\bigl(\Delta(\mu_1-\la_3,\mu_1-\la_2,
		\mu_1-\la_1),
		\Delta(\mu_1,\mu_1-\mu_2)\bigr).
		$
	\end{lemma}
	\begin{proof}
		Consider the permutation matrix $w_0 \in G$ corresponding to the longest element $(13)$ of the symmetric group $\mathfrak{S}_3$ and the automorphism $\theta$ of $G$ given by $\theta(g)=w_0(g^{-1})^T w_0$, where $g \in G$ and $(g^{-1})^T $ denotes the transpose of $g^{-1}$. For a $G$-module $M$, let ${}^\theta M$ denote the twist of $M$ by the map $\theta$. From \cite[II.2.13, after equation (1)]{Ja}, we have ${}^\theta\D(\nu)\cong\D(-\nu_3,-\nu_2,-\nu_1)$ for every dominant weight $\nu$. Thus, 	\begin{align*}
			\Hom_G(\Delta(\la),\Delta(\mu))
			&\cong
			\Hom_G\bigl(
			\Delta(-\la_3,-\la_2,-\la_1),
			\Delta(0,-\mu_2,-\mu_1)\bigr)\\
			&\cong
			\Hom_G\bigl(
			\Delta(\mu_1-\la_3,\mu_1-\la_2,
			\mu_1-\la_1),
			\Delta(\mu_1,\mu_1-\mu_2)\bigr),
		\end{align*} where the last isomorphism comes from tensoring both modules with ${\det }^{\mu_1}$.
	\end{proof}
In the previous lemma, both $(\mu_1-\la_3,\mu_1-\la_2,
	\mu_1-\la_1)$ and $(\mu_1,\mu_1-\mu_2)$ are partitions and we have $\mu_1-\mu_2 \le \mu_1-\la_3$, since $\la_3\le\la_1<\mu_2$. Hence, when studying $\Hom_G(\D(\la), \D(\mu))$, we may assume $\la$ and $\mu$ are partitions such that \begin{equation}\label{assumptions1} |\la|=|\mu|, \ \la=(\la_1,\la_2,\la_3), \ \mu=(\mu_1, \mu_2), \ \mu_2 \le \la_1,\end{equation} where $|\nu|$ denotes the sum of the parts of a partition $\nu$. We use throughout the notation \begin{equation}\label{assumptions2}\ell=\max(0,\mu_2-\la_3), \ L=\min(\la_2,\mu_2).\end{equation}

	\subsection{Basis of $\Hom_G(D(\la),\D(\mu))$}\label{secbasis}We assume throughout \eqref{assumptions1} and \eqref{assumptions2}.

	Let $E=\K^3$ be the natural $\GL_3(\K)$-module. The divided power algebra $DE = \bigoplus_{i\ge0}D_iE$ of $E$ is defined
as the graded dual of the bialgebra $S(E^*)$, where $E^*$ is
the linear dual of $E$ and $S(E^*)$ is the symmetric algebra
of $E^*$, see \cite[Section 1.1]{W}. For a sequence $\alpha=(\alpha_1, \alpha_2, \alpha_3)$ of nonnegative integers, let $D(\alpha)= D_{\alpha_1}E \otimes D_{\alpha_2}E \otimes D_{\alpha_3}E$. As a $G$-module, $D(\alpha)$ is cyclic, generated by $e_1^{(\alpha_1)}\otimes e_2^{(\alpha_2)} \otimes e_3^{(\alpha_3)}$, where $\{e_1,e_2,e_3)$ is the usual basis of $E$ \cite[isomorphism (8.3.17)]{Kr}.

For a partition $\la$, we recall that there is a natural surjection $\pi_\la:D(\la)\to\Delta(\la)$, unique up to a nonzero scalar multiple \cite[Proposition 2.1.15(a)]{W}. A map $\phi\in\Hom_G(D(\la),\Delta(\mu))$ descends to a map $\overline{\phi}\in\Hom_G(\Delta(\la),\Delta(\mu))$ if and only if $\phi(\ker\pi_\la)=0$. In this case, $\overline{\phi}$ is uniquely determined by $\phi=\overline{\phi}\circ\pi_\la$.

We use the ordered alphabet $1<2<3$. For a partition $\mu=(\mu_1,\mu_2)$, a tableau of shape $\mu$ is a filling of the diagram of $\mu$ with entries from $\{1,2,3\}$. Such a tableau is called semistandard if the entries are weakly increasing in the rows from left to right and strictly increasing in the columns from top to bottom. The weight of a tableau $T$ is the tuple $\alpha=(\alpha_1,\alpha_2,\alpha_3)$, where $\alpha_i$ is the number of appearances of $i$ in $T$. To a tableau $T$ of shape $\mu=(\mu_1, \mu_2)$ we may associate an element $x_T=x_T(1) \otimes x_T(2) \in D(\mu_1,\mu_2)$, where $x_T(i)=e_1^{(a_{i1})}e_2^{(a_{i2})}e_3^{(a_{i3})}$ and $a_{ij}$ is equal to the number of appearances of $j$ in the $i$th row of $T$, $1\le i\le2$, $1\le j\le3$.

A basis of $\D(\mu)$ is given by the elements $\pi_{\mu}(x_T)$, where $T$ runs over the semistandard tableaux of shape $\mu$ \cite[Proposition 2.1.15(b)]{W}.

A basis of $\Hom_G(D(\la),\D(\mu))$ is parameterized by the semistandard tableaux of shape $\mu$ that have weight $\la$ \cite[Lemma 8.3.18 and Corollary 8.4.14]{Kr}.
For such a tableau, let $u$ denote the number of entries equal to $2$ in its second row.
\begin{lemma}\label{basis} If $\ell > L$, then $\Hom_G(D(\la),\D(\mu))=0$. If $\ell \le L$, then a basis of $\Hom_G(D(\la),\D(\mu))$ is $\{\phi_\ell, \dots, \phi_L\}$, where for each $\ell \le u \le L$, the map $\phi_u: D(\la) \to \D(\mu)$ is defined by \[\phi_u(e_1^{(\la_1)}\otimes e_2^{(\la_2)} \otimes e_3^{(\la_3)})=\pi_\mu(e^{(\la_1)}_1e^{(\la_2-u)}_2e^{(\la_3-\mu_2+u)}_3 \otimes e^{(u)}_2e^{(\mu_2-u)}_3).\]
\end{lemma}
\begin{proof}
	This follows from \cite[Lemma 2.3 and Remark 3.3(1)]{MS} for $m=3$ and $n=3$.
\end{proof}

\subsection{The matrices $A$ and $B$}We assume throughout \eqref{assumptions1} and \eqref{assumptions2}.

Suppose $\ell \le L$. In this subsection, we define integer matrices $A$ and $B$ whose reductions modulo $2$ determine 	$\Hom_G\bigl(\D(\la),\D(\mu)\bigr)$. The columns of both are indexed by $u$ satisfying $\ell\le u\le L$. The rows are indexed by the ordered sets
\begin{align}\label{RA}
	R_A&=\{(t,k):1\le t\le\la_2,\ \ell\le k\le\min(\la_2-t,\mu_2)),\\\label{RB}
	R_B&=\{(t,k):1\le t\le\la_3,\ \max(0,\mu_2-\la_3+t)\le k\le\min(\la_2+t,\mu_2)),
\end{align}
where we consider the lexicographic ordering on pairs $(t,k)$.

For binomial coefficients we use the convention $\tbinom{a}{b}=0$ if $b<0$ or $a<b$.

\begin{definition}\label{defAB} Assume $\ell \le L$. For $\ell\le u\le L$ define
	\begin{align*}
		A_{(t,k),u}
		&=(-1)^{u-k}\tbinom{\mu_2-k}{u-k}
		\tbinom{\la_1+t-u}{t+k-u}, \ (t,k) \in R_A,\\
		B_{(t,k),u}
		&=\tbinom{k}{u}\tbinom{\la_2+t-k}{t+u-k}, \ (t,k) \in R_B.
	\end{align*}\end{definition}

We denote the reductions mod 2 of $A$ and $B$ again by $A$ and $B$. Let $\left[\begin{smallmatrix}A\\B\end{smallmatrix}\right]$ be the matrix whose rows are those of $A$ and $B$, and let $\ker_{\mathbb K}\left[\begin{smallmatrix}A\\B\end{smallmatrix}\right]$ denote its kernel over $\K$.
In the next proposition we use the notation of Section \ref{secbasis}. \begin{proposition}\label{keriso} The map $
	\ker_{\K}\left[\begin{smallmatrix}A\\B\end{smallmatrix}\right]
	\to
	\Hom_G\bigl(\D(\la),\D(\mu)\bigr)
	$ is an isomorphism, where
\[
c=(c_\ell,\ldots,c_L)^T\mapsto\overline{\phi}_c,
\ \
\overline{\phi}_c\circ\pi_\la
=\sum_{u=\ell}^{L}c_u\phi_u.
\]\end{proposition}

\begin{proof}Let
$\phi=\sum_{u=\ell}^{L}c_u\phi_{u} \in \Hom_G(D(\la),\D(\mu))
$, $c_u \in \K$.
From the second conclusion of \cite[Theorem II.3.16]{ABW} and the fact that $D(\alpha)$ is a cyclic $G$-module for any sequence $\alpha=(\alpha_1, \alpha_2, \alpha_3)$ of nonnegative integers, it follows that the kernel of the natural projection $\pi_\la: D(\la) \to \D(\la)$ is generated as a $G$-module by
\[
\begin{aligned}
	x_t&=e_1^{(\la_1)}\otimes
	e_1^{(t)}e_2^{(\la_2-t)}\otimes e_3^{(\la_3)},
	&&1\le t\le\la_2,\\
	y_t&=e_1^{(\la_1)}\otimes
	e_2^{(\la_2)}\otimes e_2^{(t)}e_3^{(\la_3-t)},
	&&1\le t\le\la_3.
\end{aligned}
\]
Thus $\phi$ descends to a map defined on \(\Delta(\la)\) if and only if $\phi(x_t)=0$ for $1\le t\le\la_2$ and $\phi(y_t)=0$ for $1\le t\le\la_3$.

If $(t,k)\in R_A$, let $T^A_{t,k}$ be the semistandard tableau of shape $\mu$ and weight $(\la_1+t,\la_2-t,\la_3)$ with second row $2^{(k)}3^{(\mu_2-k)}$. Likewise, if $(t,k)\in R_B$, let $T^B_{t,k}$ be the semistandard tableau of shape $\mu$ and weight $(\la_1,\la_2+t,\la_3-t)$ with second row $2^{(k)}3^{(\mu_2-k)}$.  The last formula for $\varphi_T(x)$ on p. 279
and the formula for $\varphi_T(y)$ immediately
preceding (2) on p. 280 of \cite{MS},
specialized to the case $m=3$ parts, give
\[
	\phi(x_t)
	=\sum_{\substack{k\\(t,k)\in R_A}}
	(Ac)_{(t,k)}\pi_\mu(x_{T^A_{t,k}}),\ \ \
	\phi(y_t)
	=\sum_{\substack{k\\(t,k)\in R_B}}
	(Bc)_{(t,k)}\pi_\mu(x_{T^B_{t,k}}).
\]
Therefore, the map $\phi: D(\la) \to \D(\mu)$ descends to a map  $\D(\la) \to \D(\mu)$ if and only if $Ac=Bc=0$. Since the maps $\phi_{u}$, $\ell \le u \le L$, form a basis of $\Hom_G(D(\la),\D(\mu))$, the result follows. \end{proof}

	\subsection{Binomial congruences}\label{secbinom}

	We will use the following special case of Lucas' theorem.
	\begin{lemma}\label{lembinom1}
		Let $N\ge0$, $M\in\mathbb Z$ and
		$\epsilon,\eta\in\{0,1\}$. Then
	$
		\tbinom{2N+\epsilon}{2M+\eta}
		\equiv\tbinom NM\tbinom\epsilon\eta\pmod2.$
	\end{lemma}

	For a positive integer $n$, let $v_2(n)$ be the exponent of $2$ in the prime factorization of
	$n$.

	\begin{lemma}\cite[Corollary 22.5]{Jam}\label{lembinom2}
Let $n,d$ be nonnegative integers. Then
$
\tbinom{n+t}{t}\equiv0\pmod2
\ \text{for every }1\le t\le d
$
if and only if $d<2^{v_2(n+1)}$.
	\end{lemma}

\section{The recursion}\label{sec3}
Throughout this section we assume \eqref{assumptions1} and \eqref{assumptions2}. Also we assume $\ell \le L$.

\subsection{Reduction to $(0,0,0)$ and $(0,1,1)$}Suppose $\ell<L$.
Let $c=(c_\ell, \dots, c_L) \in \ker_{\mathbb K}\left[\begin{smallmatrix}A\\B\end{smallmatrix}\right]$. First we establish that $(1,u)\in R_A$ and $(1,u+1)\in R_B$ for every $\ell\le u<L$. Since $\ell<L$, we have $\la_3>0$. For
$\ell\le u<L$, we have
$
\ell\le u\le L-1\le\min(\la_2-1,\mu_2)
$
and
\[
\max(0,\mu_2-\la_3+1)
\le\ell+1\le u+1
\le L\le\min(\la_2+1,\mu_2).
\]
Thus we obtain $(1,u)\in R_A$ and $(1,u+1)\in R_B$
by \eqref{RA} and \eqref{RB}, respectively. 

For $\ell \le u < L$, the rows $(1,u)$ of $A$
and $(1,u+1)$ of $B$ give the following equations
	\begin{align}\label{eq31}
		(\la_1+1-u)c_u+(\mu_2-u)c_{u+1}&=0,\\\label{eq32}
		(u+1)c_u+(\la_2-u)c_{u+1}&=0.
	\end{align}
Adding these and using that the characteristic is $2$, we have \begin{equation}\label{eq33}\la_1c_u+(\la_2+\mu_2)c_{u+1}=0.\end{equation}
In the next lemma, we consider all possibilities modulo 2 of the triple $(\la_1,\la_2,\mu_2)$.
\begin{lemma}\label{lemmod2}
	Suppose $\ell<L$ and $Ac=Bc=0$, where $c=(c_\ell, \dots, c_L)$. Then $c$ satisfies the following conditions, according to the parity of $(\la_1,\la_2,\mu_2)$.
	\[
	\begin{array}{c|l}
		(\la_1,\la_2,\mu_2)\bmod2
		&\text{Consequence for }c\\\hline
		(0,0,0)
		&c_u=0\quad(\ell\le u\le L,\ u\text{ even}),\\
		(0,0,1),(0,1,0)
		&c_u=0\quad(\ell<u\le L),\\
		(0,1,1)
		&c_u=c_{u+1}\quad(\ell\le u<L,\ u\text{ even}),\\
		(1,0,0),(1,0,1),(1,1,1)
		&c_u=0\quad(\ell\le u\le L),\\
		(1,1,0)
		&c_u=c_{u+1}\quad(\ell\le u<L).
	\end{array}
	\]
	Thus, $\dim_\K\Hom_G(\D(\la),\D(\mu)) \le 1$ in all cases except possibly when $(\la_1, \la_2, \mu_2) \equiv (0,0,0) \ \text{or} \ (0,1,1) \mod 2$.
\end{lemma}
\begin{proof}
	For the triples $(1,0,0)$ and $(1,1,1)$ in the table, we have $c=0$.
	Indeed, \eqref{eq33} gives $c_u=0$ for $\ell\le u<L$.
	Since $L=\min(\la_2,\mu_2)$ and
	$\la_2\equiv\mu_2\pmod2$, we have
	$\la_2-L+1\equiv1\pmod2$. Thus eq. \eqref{eq32},
	with $u=L-1$, gives $c_L=0$.
	
	For the triple $(1,0,1)$, eq. \eqref{eq33} gives
	$c_u=c_{u+1}$ for $\ell\le u<L$.
	Substituting in eq. \eqref{eq31} gives $c_u=0$
	for $\ell\le u<L$, and hence $c=0$.
	
	For the triple $(0,0,0)$, if $L$ is even,
	equation \eqref{eq31} with $u=L-1$ gives $c_L=0$.
	
	The remaining assertions in the table follow from
	\eqref{eq31}-\eqref{eq33}.
	
	The final conclusion of the lemma follows from 
	Proposition \ref{keriso}.
	\end{proof}

\subsection{Halving procedure}	According to the previous lemma, we have $\dim_\K\Hom_G(\D(\la),\D(\mu)) \le 1$ in all cases except possibly  the cases $(0,0,0)$ and $(0,1,1)$ for $(\la_1, \la_2, \mu_2) \mod 2$. For these we will apply a recursion that relates $\Hom_G(\D(\la),\D(\mu))$ to $\Hom_G(\D(\la'),\D(\mu'))$ for certain smaller partitions $\la'$, $\mu'$.

Let us give a roadmap of the present subsection. For each pair of partitions $(\la, \mu)$ that belong to the above cases with $\ell < L$, we define a new pair of partitions $(\la', \mu')$ (Definition \ref{defhalf}). Let $A', B'$ be their corresponding matrices of Definition \ref{defAB}. In Lemma \ref{lemsmaller} we show that $(\la',\mu')$ again satisfies \eqref{assumptions1}, that $\ell'\le L'$ and that $\la'_3<\la_3$. In Lemma \ref{lemrowhalf} we relate certain row indices of the matrices $A, B$ with those of $A', B'$, to be used in the proof of Proposition \ref{matrixreduction}. In that proposition we apply column operations to $A, B$ and, using Lucas' theorem, obtain a pair of matrices $\widetilde A,\widetilde B$ such that $\widetilde A$ and $A'$ have the same nonzero rows (without counting repetitions) and likewise for $\widetilde B$ and $B'$. Finally, in Proposition \ref{propiso} we use Lemma \ref{lemmod2} and the matrices $\widetilde A, \widetilde B$ to establish an isomorphism $\ker_\K\left[\begin{smallmatrix}A\\B\end{smallmatrix}\right]
\cong \ker_\K\left[\begin{smallmatrix}A'\\B'\end{smallmatrix}\right]$, which by Proposition \ref{keriso} gives the desired isomorphism $\Hom_G(\D(\la),\D(\mu)) \cong \Hom_G(\D(\la'),\D(\mu'))$.

For a rational number $a$, let $\lfloor a \rfloor$ be the largest integer $m$ such that $m \le a$.
	\begin{definition}\label{defhalf}
		Suppose $\ell<L$, $\la_1$ is even, and
		$\la_2\equiv\mu_2\equiv\varepsilon\pmod2$, where
		$\varepsilon\in\{0,1\}$. Define
		\begin{align*}
			\la'_1&=\la_1/2-1,&
			\la'_2&=\lfloor(\la_2-1)/2\rfloor,&
			\la'_3&=\lfloor(\la_3-1+\varepsilon)/2\rfloor,\\
			\mu'_2&=\lfloor(\mu_2-1)/2\rfloor,&
			\mu'_1&=\la'_1+\la'_2+\la'_3-\mu'_2.
		\end{align*}
		We say that $\la'=(\la'_1,\la'_2,\la'_3)$ and
		$\mu'=(\mu'_1,\mu'_2)$ are obtained from $\la,\mu$ by
		\emph{halving}.

		Also define $\ell'=\max(0,\mu'_2-\la'_3)$ and
		$L'=\min(\la'_2,\mu'_2)$.
	\end{definition}
	
	To maintain the flow of the argument, proofs of Lemma \ref{lemsmaller} and Lemma \ref{lemrowhalf} below are deferred to Section \ref{lemproofs}.
	\begin{lemma}\label{lemsmaller}Under the assumptions of Definition \ref{defhalf}, the following hold.
		\begin{enumerate}
			\item $\la'$ and $\mu'$ are partitions with $|\la'|=|\mu'|$. Also
			$0\le\mu'_2\le\la'_1$ and
			$0\le\la'_3<\la_3$.
			\item $\ell'=\lfloor\ell/2\rfloor$ and
			$L'=\lfloor(L-1)/2\rfloor$. In particular, $\ell'\le L'$.
			\item We have $\ell\le2i+1\le L$ if and only if $\ell'\le i\le L'$.
		\end{enumerate}
	\end{lemma}

For the partitions $\la',\mu'$ of Lemma \ref{lemsmaller}, let $A'$ and $B'$ be the corresponding matrices of Definition \ref{defAB} and let $R_{A'}$ and $R_{B'}$ be the corresponding sets of row indices given by \eqref{RA} and \eqref{RB}.

\begin{lemma}\label{lemrowhalf}
Under the assumptions of Definition \ref{defhalf}, the following hold. Let $h\ge1$ and $j \ge 0$.
\begin{enumerate}
\item
$
(2h,2j+1)\in R_A \ \Leftrightarrow\ (h,j)\in R_{A'}, \ \ \text{and} \ \
(2h,2j+1)\in R_B \ \Leftrightarrow\ (h,j)\in R_{B'}.
$
\item Let $\varepsilon=1$. Then
\[
\begin{aligned}
(2h,2j)\in R_A&\ \Rightarrow\ (2h,2j+1)\in R_A,\\
(2h,2j)\in R_B&\ \Rightarrow\ (2h,2j+1)\in R_B.
\end{aligned}
\]
In particular, if $(2h,2j)\in R_A$, then $(h,j)\in R_{A'}$; if
$(2h,2j)\in R_B$, then $(h,j)\in R_{B'}$.
\end{enumerate}
\end{lemma}

\begin{proposition}\label{matrixreduction}
	Under the assumptions of Definition \ref{defhalf}, perform the
	following operations on each of $A, B$.
	\begin{enumerate}
		\item If $\varepsilon=0$, delete all columns with even index.
		\item If $\varepsilon=1$, add every column with even index to the
		column on the right, then delete the columns with even indices.
		An odd-indexed column is left unchanged if the preceding
		even-indexed column is absent.
	\end{enumerate}
Reindex the remaining columns so that the column originally indexed by \(2i+1\) is now indexed by \(i\). Denote the resulting matrices by \(\widetilde A,\widetilde B\). Then every row of $\widetilde A$ is zero
	or a row of $A'$, and every row of $A'$ occurs in $\widetilde A$. The same holds
	for $\widetilde B$ and $B'$.
\end{proposition}

\begin{proof}
	(1) Suppose $\varepsilon=0$. 
	
	\textit{Case 1A}. We show that every nonzero row of
	$\widetilde A$ is a row of $A'$.
	Consider such a row, indexed by $(t,k)\in R_A$.
	Throughout the proof, $i$ ranges over all integers with
	$\ell'\le i\le L'$, which index the columns of $\widetilde A$
	and $\widetilde B$ by Lemma \ref{lemsmaller}(3).
	Recall from Definition \ref{defhalf} that
	\begin{equation}\label{lala'}
	\la_1=2\la'_1+2,\qquad
	\la_2=2\la'_2+2,\qquad
	\mu_2=2\mu'_2+2.
	\end{equation}
	From the definition of $\widetilde A$ and Definition \ref{defAB} we have
	\[
	\widetilde A_{(t,k),i}
	=A_{(t,k),2i+1}
	=\tbinom{\mu_2-k}{2i+1-k}
	\tbinom{\la_1+t-2i-1}{t+k-2i-1}.
	\]
	If $k$ is even, Lemma \ref{lembinom1} gives
$
	\tbinom{\mu_2-k}{2i+1-k}=0,
$
	so the row is zero. Hence $k$ is odd.
	If $t$ is odd, Lemma \ref{lembinom1} gives
$
	\tbinom{\la_1+t-2i-1}{t+k-2i-1}=0,
$
	so the row is again zero. So we may assume $t=2h$ and $k=2j+1$, where
	 $h\ge1$ and $j\ge0$.
	
	Since $(2h,2j+1)\in R_A$, from 
	Lemma \ref{lemrowhalf}(1) we get $(h,j)\in R_{A'}$.
	By Lemma \ref{lembinom1} and \eqref{lala'},
	\begin{equation}\label{rowAepsilonzero}
		\
			\widetilde A_{(2h,2j+1),i}
			=\tbinom{\mu_2-2j-1}{2i-2j}
			\tbinom{\la_1+2h-2i-1}{2h+2j-2i}=\tbinom{\mu'_2-j}{i-j}
			\tbinom{\la'_1+h-i}{h+j-i}
			=A'_{(h,j),i}.	\end{equation}
	This proves Case 1A.
	
	\textit{Case 1B}. We show that every nonzero row of
	$\widetilde B$ is a row of $B'$.
	Consider such a row, indexed by $(t,k)\in R_B$. We have
	\[
	\widetilde B_{(t,k),i}
	=B_{(t,k),2i+1}
	=\tbinom{k}{2i+1}
	\tbinom{\la_2+t-k}{t+2i+1-k}.
	\]
	If $k$ is even, Lemma \ref{lembinom1} gives
$
	\tbinom{k}{2i+1}=0,
$
	so the row is zero. Hence $k$ is odd.
	If $t$ is odd, Lemma \ref{lembinom1} gives
$
	\tbinom{\la_2+t-k}{t+2i+1-k}=0,
	$
	so the row is again zero. Thus $t=2h$ and $k=2j+1$, where $h\ge1$ and $j\ge0$.
	
	Since $(2h,2j+1)\in R_B$ from
	Lemma \ref{lemrowhalf}(1) we get $(h,j)\in R_{B'}$.
	Again from Lemma \ref{lembinom1} we get
	\begin{equation}\label{rowBepsilonzero}
			\widetilde B_{(2h,2j+1),i}
			=\tbinom{2j+1}{2i+1}
			\tbinom{\la_2+2h-2j-1}{2h+2i-2j}=\tbinom ji\tbinom{\la'_2+h-j}{h+i-j}
			=B'_{(h,j),i}.
	\end{equation}
	This proves Case 1B.
	
	\textit{Conclusion of (1)}. Every row of $A'$ occurs in
	$\widetilde A$, and every row of $B'$ occurs in $\widetilde B$.
	Indeed, let $(h,j)\in R_{A'}$. By Lemma \ref{lemrowhalf}(1) we have
	$(2h,2j+1)\in R_A$, and \eqref{rowAepsilonzero} gives
	$\widetilde A_{(2h,2j+1),i}=A'_{(h,j),i}$.
	If $(h,j)\in R_{B'}$, the same lemma gives
	$(2h,2j+1)\in R_B$, and \eqref{rowBepsilonzero} gives
	$\widetilde B_{(2h,2j+1),i}=B'_{(h,j),i}$.
	
	(2) Suppose $\varepsilon=1$. 	From Definition \ref{defhalf} we have
	\begin{equation}\label{lala'2}
	\la_1=2\la'_1+2,\qquad
	\la_2=2\la'_2+1,\qquad
	\mu_2=2\mu'_2+1.
	\end{equation}
	
	\textit{Case 2A}. We show that every nonzero row of
	$\widetilde A$ is a row of $A'$.
	Consider such a row, indexed by $(t,k)\in R_A$.
	Since $L=\min(\la_2,\mu_2)$ is odd, every even-indexed
	column of $A$ or $B$ has its following odd-indexed column.
	
	If column $2i$ is absent, then $2i=\ell-1$ and $\ell>0$,
	since the original column indices are consecutive.
	By \eqref{assumptions2}, $\ell=\mu_2-\la_3$.
	Since $(t,k)\in R_A$, equation \eqref{RA} gives $k\ge\ell$, so
	\[
	2i-k=\ell-1-k<0,\qquad
	\tbinom{\mu_2-k}{2i-k}=0.
	\]
	Thus, whether or not column $2i$ is present, Definition \ref{defAB}
	gives
	\begin{equation}\label{rowAtilde}
		\widetilde A_{(t,k),i}
		=\tbinom{\mu_2-k}{2i-k}
		\tbinom{\la_1+t-2i}{t+k-2i}+\tbinom{\mu_2-k}{2i+1-k}
		\tbinom{\la_1+t-2i-1}{t+k-2i-1}.
	\end{equation}
	
	Write $k=2j$ or $k=2j+1$, where $j\ge0$ is an integer.
	If $k=2j$, Lemma \ref{lembinom1} gives
	\[
	\tbinom{\mu_2-2j}{2i-2j}
	=\tbinom{\mu_2-2j}{2i+1-2j}
	=\tbinom{\mu'_2-j}{i-j}.
	\]
	From Pascal's identity
	in characteristic $2$ and \eqref{rowAtilde} we have
	\[
		\widetilde A_{(t,2j),i}
		=\tbinom{\mu'_2-j}{i-j}
		\left[
		\tbinom{\la_1+t-2i}{t+2j-2i}
		+\tbinom{\la_1+t-2i-1}{t+2j-2i-1}
		\right]=\tbinom{\mu'_2-j}{i-j}
		\tbinom{\la_1+t-2i-1}{t+2j-2i}.
	\]
	If $k=2j+1$, Lemma \ref{lembinom1} and a similar computation give
	\[
	\begin{gathered}
		\tbinom{\mu_2-2j-1}{2i-2j-1}=0,\ \
		\tbinom{\mu_2-2j-1}{2i-2j}=\tbinom{\mu'_2-j}{i-j},\ \
		\widetilde A_{(t,2j+1),i}
		=\tbinom{\mu'_2-j}{i-j}
		\tbinom{\la_1+t-2i-1}{t+2j-2i}.
	\end{gathered}
	\]
	Now we conclude that for both cases of $k$, if $t$ is odd, Lemma \ref{lembinom1}
	gives
$
	\tbinom{\la_1+t-2i-1}{t+2j-2i}=0
$
	anb thus the row is zero. Thus we have $t=2h$, where $h\ge1$.
	
	If $k\in\{2j,2j+1\}$ with
	$(2h,k)\in R_A$, then we have $(h,j)\in R_{A'}$ by
	Lemma \ref{lemrowhalf}(1) when $k=2j+1$, and by
	Lemma \ref{lemrowhalf}(2) when $k=2j$.
	Since $\la_1+2h-2i-1=2(\la'_1+h-i)+1$ is odd,
	Lemma \ref{lembinom1} gives
	\begin{equation}\label{rowAepsilonone}
			\widetilde A_{(2h,k),i}
			=\tbinom{\mu'_2-j}{i-j}
			\tbinom{\la_1+2h-2i-1}{2h+2j-2i}=\tbinom{\mu'_2-j}{i-j}
			\tbinom{\la'_1+h-i}{h+j-i}
			=A'_{(h,j),i}
	\end{equation}
	and we have proven the claim of Case 2A.
	
	\textit{Case 2B}. We show that every nonzero row of
	$\widetilde B$ is a row of $B'$.
	Consider such a row, indexed by $(t,k)\in R_B$, and write
	$k=2j$ or $k=2j+1$, where $j\ge0$ is an integer.
	If column $2i$ is absent, then $2i=\ell-1$ and
	$\ell=\mu_2-\la_3>0$, as above.
	Equation \eqref{RB} gives $k\ge\ell+t$, so
	\[
	t+2i-k=t+\ell-1-k<0,\qquad
	\tbinom{\la_2+t-k}{t+2i-k}=0.
	\]
	Thus Definition \ref{defAB} gives, whether or not column $2i$
	is present,
		\begin{equation}\label{rowBtilde}
		\widetilde B_{(t,k),i}
		=\tbinom{k}{2i}
		\tbinom{\la_2+t-k}{t+2i-k}+\tbinom{k}{2i+1}
		\tbinom{\la_2+t-k}{t+2i+1-k}.
	\end{equation}
	If $k=2j$, Lemma \ref{lembinom1} gives
	\[
	\begin{gathered}
		\tbinom{2j}{2i+1}=0,\ \tbinom{2j}{2i}=\tbinom ji,\ \
		\widetilde B_{(t,2j),i}
		=\tbinom ji\tbinom{\la_2+t-2j}{t+2i-2j}.
	\end{gathered}
	\]
	If $k=2j+1$, Lemma \ref{lembinom1} gives
$
	\tbinom{2j+1}{2i}
	=\tbinom{2j+1}{2i+1}
	=\tbinom ji
$
	and from
	Pascal's identity and \eqref{rowBtilde} we get
	\[
		\widetilde B_{(t,2j+1),i}
		=\tbinom ji\left[
		\tbinom{\la_2+t-2j-1}{t+2i-2j-1}
		+\tbinom{\la_2+t-2j-1}{t+2i-2j}
		\right]=\tbinom ji\tbinom{\la_2+t-2j}{t+2i-2j}.
	\]
	For either value of $k$, if $t$ is odd, Lemma \ref{lembinom1}
	gives
$
	\tbinom{\la_2+t-2j}{t+2i-2j}=0,
$
	so the row is zero. Thus we have $t=2h$, where $h\ge1$.
	
	If $k\in\{2j,2j+1\}$ with
	$(2h,k)\in R_B$, then we have $(h,j)\in R_{B'}$ by
	Lemma \ref{lemrowhalf}(1) when $k=2j+1$, and by
	Lemma \ref{lemrowhalf}(2) when $k=2j$.
	Since $\la_2+2h-2j=2(\la'_2+h-j)+1$ is odd,
	Lemma \ref{lembinom1} gives
	\begin{equation}\label{rowBepsilonone}
			\widetilde B_{(2h,k),i}
			=\tbinom ji\tbinom{\la_2+2h-2j}{2h+2i-2j}=\tbinom ji\tbinom{\la'_2+h-j}{h+i-j}
			=B'_{(h,j),i}
	\end{equation}
and we have proven the claim of Case 2B.

	\textit{Conclusion of (2)}. Every row of $A'$ occurs in
	$\widetilde A$, and every row of $B'$ occurs in $\widetilde B$.
	Let $(h,j)\in R_{A'}$. By Lemma \ref{lemrowhalf}(1) we have
	$(2h,2j+1)\in R_A$, and \eqref{rowAepsilonone} gives
	$\widetilde A_{(2h,2j+1),i}=A'_{(h,j),i}$.
	If $(h,j)\in R_{B'}$, the same lemma gives
	$(2h,2j+1)\in R_B$, and from \eqref{rowBepsilonone} we get
	$\widetilde B_{(2h,2j+1),i}=B'_{(h,j),i}$.
\end{proof}

\begin{example}
	Let $\la=(10,7,2)$ and $\mu=(14,5)$. Then $\ell=3$, $L=5$,
	and $\varepsilon=1$. For the corresponding matrices (Definition \ref{defAB}) we have
	\begin{align*}
		R_A&=\{(1,3),(1,4),(1,5),(2,3),(2,4),(2,5),(3,3),(3,4),(4,3)\},\\
		R_B&=\{(1,4),(1,5),(2,5)\}.
	\end{align*}
	The columns are indexed by $3,4,5$. Thus
	\[
	A=\begin{pmatrix}
		\binom{8}{1}&\binom{2}{1}&0\\[3pt]
		0&\binom{7}{1}&\binom{1}{1}\\[3pt]
		0&0&\binom{6}{1}\\[3pt]
		\binom{9}{2}&\binom{2}{1}\binom{8}{1}&\binom{2}{2}\\[3pt]
		0&\binom{8}{2}&\binom{1}{1}\binom{7}{1}\\[3pt]
		0&0&\binom{7}{2}\\[3pt]
		\binom{10}{3}&\binom{2}{1}\binom{9}{2}&\binom{2}{2}\binom{8}{1}\\[3pt]
		0&\binom{9}{3}&\binom{1}{1}\binom{8}{2}\\[3pt]
		\binom{11}{4}&\binom{2}{1}\binom{10}{3}&\binom{2}{2}\binom{9}{2}
	\end{pmatrix},
	\qquad
	B=\begin{pmatrix}
		\binom{4}{3}&\binom{4}{4}\binom{4}{1}&0\\[3pt]
		0&\binom{5}{4}&\binom{5}{5}\binom{3}{1}\\[3pt]
		\binom{5}{3}&\binom{5}{4}\binom{4}{1}&\binom{5}{5}\binom{4}{2}
	\end{pmatrix}.
	\]
	According to Proposition \ref{matrixreduction}(2), in each matrix, replace the column indexed by $5$ with the sum
	of the columns indexed by $4$ and $5$, and delete the column
	indexed by $4$. Reindexing the two remaining columns by
	$3\mapsto 1$ and $5\mapsto 2$ gives
	\[
	\widetilde A=\begin{pmatrix}
		0&0\\
		0&0\\
		0&0\\
		0&1\\
		0&1\\
		0&1\\
		0&0\\
		0&0\\
		0&0
	\end{pmatrix},
	\qquad
	\widetilde B=\begin{pmatrix}
		0&0\\
		0&0\\
		0&0
	\end{pmatrix}.
	\]
	On the other hand, Definition \ref{defhalf} gives the partitions $\la'=(4,3,1)$ and $\mu'=(6,2)$. Their corresponding matrices (Definition \ref{defAB}) are
	\[
	A'=\begin{pmatrix}
		\binom{4}{1}&\binom{1}{1}\\[3pt]
		0&\binom{3}{1}\\[3pt]
		\binom{5}{2}&\binom{1}{1}\binom{4}{1}
	\end{pmatrix}
	=\begin{pmatrix}
		0&1\\
		0&1\\
		0&0
	\end{pmatrix},
	\qquad
	B'=\begin{pmatrix}
		\binom{2}{1}&\binom{2}{2}\binom{2}{1}
	\end{pmatrix}
	=\begin{pmatrix}0&0\end{pmatrix}.
	\]
By Case 2A and Case 2B in the proof of Proposition \ref{matrixreduction},
row $(2,3)$ of $\widetilde A$ corresponds to row $(1,1)$ of $A'$,
rows $(2,4)$ and $(2,5)$ correspond to row $(1,2)$, and
row $(4,3)$ corresponds to row $(2,1)$.
Similarly, row $(2,5)$ of $\widetilde B$ corresponds to
row $(1,2)$ of $B'$.
\end{example}

	\begin{proposition}\label{propiso}
		Under the assumptions of Definition \ref{defhalf}, the map
		\[
		\Phi:\ker_\K\left[\begin{smallmatrix}A\\B\end{smallmatrix}\right]
		\longrightarrow\ker_\K\left[\begin{smallmatrix}A'\\B'\end{smallmatrix}\right],
		\qquad
		(c_u)_{\ell\le u\le L}\mapsto
		(c_{2i+1})_{\ell'\le i\le L'}
		\]
		is an isomorphism. Its inverse sends $c'$ to $c=(c_u)_{\ell\le u\le L}$ given by
		\begin{equation}\label{inv}
			c_u=
			\begin{cases}
				c'_{(u-1)/2},&u\text{ odd},\\
				0,&u\text{ even and }\varepsilon=0,\\
				c'_{u/2},&u\text{ even and }\varepsilon=1.
			\end{cases}
		\end{equation}
	\end{proposition}

\begin{proof}
By Proposition \ref{matrixreduction} we know that
$
\ker_\K\begin{bmatrix}
	\widetilde A\\ \widetilde B
\end{bmatrix}
=
\ker_\K\begin{bmatrix}
	A'\\B'
\end{bmatrix}.
$

	(1) We show that $\Phi$ is well defined.
	Let $c=(c_u)_{\ell\le u\le L}\in
	\ker_\K\left[\begin{smallmatrix}A\\B\end{smallmatrix}\right]$. For even $u$ with $\ell\le u\le L$, Lemma \ref{lemmod2}
	gives $c_u=0$ if $\varepsilon=0$ and $c_u=c_{u+1}$
	if $\varepsilon=1$. Thus from the definition of $\widetilde A$ and $\widetilde B$ in
	Proposition \ref{matrixreduction} we have
$
	\widetilde A\Phi(c)=Ac=0,
	\
	\widetilde B\Phi(c)=Bc=0.
	$
	Hence
$
	\operatorname{im}\Phi\subseteq
	\ker_\K\left[\begin{smallmatrix}A'\\B'\end{smallmatrix}\right].
$

(2) We show that $\Phi$ is injective. Suppose $\Phi(c)=0$.
	By Lemma \ref{lemsmaller}(3), we have $c_u=0$ for every odd $u$
	with $\ell\le u\le L$.
	Lemma \ref{lemmod2} (first and third rows of the table) then gives $c=0$.

(3) We show that $\Phi$ is surjective.
	Let $c'\in\ker_\K\left[\begin{smallmatrix}A'\\B'\end{smallmatrix}\right]$.
For $\ell\le u\le L$, Lemma \ref{lemsmaller}(3) gives
$\ell'\le (u-1)/2\le L'$ if $u$ is odd.
If $\varepsilon=1$ and $u$ is even, then $L$ is odd,
so $\ell\le u+1\le L$; the same lemma gives
$\ell'\le u/2\le L'$.
Define $c=(c_u)_{\ell\le u\le L}$ by \eqref{inv}. The column operations in Proposition \ref{matrixreduction}
show that $Ac=\widetilde A c'$ and $Bc=\widetilde B c'$.
By the same proposition, $\widetilde A c'=\widetilde B c'=0$.
Hence $c\in\ker_\K\left[\begin{smallmatrix}A\\B\end{smallmatrix}\right]$.
Since $c_{2i+1}=c'_i$ by \eqref{inv}, we have $\Phi(c)=c'$.
\end{proof}

\begin{remark}\label{steinberghalving}In Definition \ref{defhalf} suppose
	 $\varepsilon=1$ and $\la_3$ is even. Then an
isomorphism
$
\Hom_G(\D(\la),\D(\mu))$ $
\cong
\Hom_G(\D(\la'),\D(\mu'))
$ can be obtained using the Steinberg module.
Indeed, let $\rho=(2,1,0)$ and
$\mathrm{St}=\Delta(\rho)$ be the Steinberg module.
Then $\la=\rho+2\la'$ and $\mu=\rho+2\mu'$.
The Steinberg tensor identity
\cite[II.3.18, 3.19]{Ja} gives
$
\D(\la)\cong
\mathrm{St}\otimes\D(\la')^{(1)},$
 $
\D(\mu)\cong
\mathrm{St}\otimes\D(\mu')^{(1)}
$, where $M^{(1)}$ denotes the first Frobenius twist of a $G$-module $M$. The functor
$M\mapsto\mathrm{St}\otimes M^{(1)}$
is fully faithful \cite[II.10.5(1)]{Ja} and thus
we obtain the isomorphism.
\end{remark}

\section{Main result}\label{sec4}
	\subsection{A recursive nonvanishing criterion}The following is the main result of the paper.
\begin{theorem}\label{thm2}Let $\la, \mu$ satisfy \eqref{assumptions1} and let $\ell, L$ be as in \eqref{assumptions2}.
	If $\ell>L$, then
	$\Hom_G(\Delta(\la),\Delta(\mu))=0$.
	If $\ell\le L$, then exactly one of the following cases applies.
	\begin{enumerate}
		\item\label{case:constant}
		Suppose $\ell=L$ or $\la_1$ is odd.
		Then $\Hom_G(\Delta(\la),\Delta(\mu)) \neq 0$ if and only if
		\begin{equation}\label{thm(1)}
			\la_2-\ell<2^{v_2(\la_1-\mu_2+1)},
			\qquad
			\la_3<2^{v_2(\la_2+1)}.
		\end{equation}
		If \eqref{thm(1)} holds, then $\dim_\K\Hom_G(\Delta(\la),\Delta(\mu))=1$ and a generator is induced by the map $\sum_{u=\ell}^{L}\phi_u : D(\la) \to \D(\mu)$.
		\item\label{case:isolated}
		Suppose $\ell<L$, $\la_1$ is even,
		and $\la_2,\mu_2$ have opposite parity.
		Then $\Hom_G(\Delta(\la),\Delta(\mu)) \neq 0$ if and only if
		\begin{equation}\label{thm(2)}
			\la_2-\ell<2^{v_2(\la_1-\ell+1)},
			\qquad
			\la_3<2^{v_2(\ell+1)}.
		\end{equation}
		If \eqref{thm(2)} holds, then $\dim_\K\Hom_G(\Delta(\la),\Delta(\mu))=1$ and a generator is induced by the map $\phi_\ell : D(\la) \to \D(\mu)$.
		\item\label{case:halving}
		Suppose $\ell<L$, $\la_1$ is even,
		and $\la_2,\mu_2$ have the same parity.
		For the partitions $\la',\mu'$ of
		Definition \ref{defhalf} we have
	$\Hom_G\bigl(\Delta(\la),\Delta(\mu)\bigr) \cong
		\Hom_G\bigl(\Delta(\la'),\Delta(\mu')\bigr).
	$ Moreover, if $\sum_{i=\ell'}^{L'}c'_i\phi'_i$ induces a generator of
	$\Hom_G\bigl(\Delta(\la'),\Delta(\mu')\bigr)$, where $\phi'_i$ are the maps of
	Lemma \ref{basis} for $(\la',\mu')$, then, with $c_u$ defined by \eqref{inv},
	$\sum_{u=\ell}^{L}c_u\phi_u$ induces a generator of
	$\Hom_G\bigl(\Delta(\la),\Delta(\mu)\bigr)$.
	\end{enumerate}
Repeated application of (3) terminates in (1) or (2). Thus, $\dim_\K\Hom_G(\Delta(\la),\Delta(\mu)) \leq 1$.
\end{theorem}

\begin{proof}
	By Lemma \ref{basis}, 	we have $\Hom_G(\Delta(\la),\Delta(\mu)) = 0$
	if $\ell>L$. Let $\ell\le L$.
	
	(1) Suppose $\ell=L$ or $\la_1$ is odd.
	If $\ell=L$, the matrix
	$\left[\begin{smallmatrix}A\\B\end{smallmatrix}\right]$
	has one column. If $\ell<L$, then $\la_1$ is odd, and
	Lemma \ref{lemmod2} shows that every vector in its kernel
	has equal coordinates (possibly all zero).
	Thus in both cases the kernel is
	contained in the subspace of  $\K^{L-\ell+1}$ spanned by
	$(1,\ldots,1)^T$.
	By Proposition \ref{keriso}, $\Hom_G(\Delta(\la),\Delta(\mu))$ is
	nonzero if and only if every row of $A$ and
	$B$ has sum zero. For $(t,k)\in R_A$, we compute
	\begin{equation}\label{rowsumA}
		\sum_{u=\ell}^{L}A_{(t,k),u}
		=\sum_{u=\ell}^{L}(-1)^{u-k}
		\tbinom{\mu_2-k}{u-k}
		\tbinom{\la_1+t-u}{t+k-u}=\sum_{j=0}^{t}(-1)^j
		\tbinom{\mu_2-k}{j}
		\tbinom{\la_1+t-k-j}{t-j}
		=\tbinom{\la_1-\mu_2+t}{t}.
	\end{equation}
	In \eqref{rowsumA} we used the change of variable $j=u-k$. So 
	$\ell-k\le j\le L-k$. For the new end points, by \eqref{RA}, we have $\ell\le k$ and
	$t\le\la_2-k$. Hence 
	\[
	\ell-k\le j\le L-k,\quad 0\le j\le t
	\iff 0\le j\le\min(t,L-k)
	\iff 0\le j\le\min(t,\mu_2-k),
	\] where in the last equivalence we used $L=\min(\lambda_2,\mu_2)$.
	Since the summand $(-1)^j
	\tbinom{\mu_2-k}{j}
	\tbinom{\la_1+t-k-j}{t-j}$ vanishes for $j<0$, $j>t$, or $j>\mu_2-k$ we may sum over $0\le j\le t$. The last equality in \eqref{rowsumA} uses
	\cite[Lemma 4.1(2)]{MS} with
	$a=\la_1-k$, $b=\mu_2-k$, and $c=t$.
	
	For $(t,k)\in R_B$, a similar computation yields
	\begin{equation}\label{rowsumB}
			\sum_{u=\ell}^{L}B_{(t,k),u}
			=\sum_{u=\ell}^{L}
			\tbinom{k}{u}
			\tbinom{\la_2+t-k}{t+u-k}=\sum_{j=0}^{t}
			\tbinom{k}{j}
			\tbinom{\la_2+t-k}{t-j}=\tbinom{\la_2+t}{t}.
	\end{equation}
Indeed, we used $j=k-u$, so
$k-L\le j\le k-\ell$.
By \eqref{RB}, we have $k\le\mu_2$ and
$k-\ell\ge\min(k,t)$.
Hence
\[
\begin{aligned}
	k-L\le j\le k-\ell,\ \ 
	0\le j\le\min(k,t)&\iff
	\max(0,k-L)\le j\le\min(k,t)\\
	&\iff
	\max(0,k-\la_2)\le j\le\min(k,t),
\end{aligned}
\]
where the last equivalence uses
$L=\min(\la_2,\mu_2)$ and $k\le\mu_2$.
The summand 	$\tbinom{k}{j}
\tbinom{\la_2+t-k}{t-j}$ vanishes for $j<0$, $j>k$,
$j>t$, or $j<k-\la_2$.
Thus we may sum over $0\le j\le t$.
The last equality in \eqref{rowsumB}
is Vandermonde's identity.
	
By \eqref{RA} and \eqref{RB}, the ranges of $t$
are $1\le t\le\la_2-\ell$ and
$1\le t\le\la_3$, respectively.
	Thus equations \eqref{rowsumA} and \eqref{rowsumB} show that every row
	of $A$ and $B$ has sum zero if and only if
	\[
	\tbinom{\la_1-\mu_2+t}{t}
	\equiv0\pmod2 \ \
	(1\le t\le\la_2-\ell),\ \ \
	\tbinom{\la_2+t}{t}
	\equiv0\pmod2 \ \
	(1\le t\le\la_3).
	\]
	By Lemma \ref{lembinom2}, these conditions
	are equivalent to \eqref{thm(1)}.
	When they hold, the vector $(1,\ldots,1)^T$ spans 	$\ker_\K\left[\begin{smallmatrix}A\\B\end{smallmatrix}\right]$. Then Proposition \ref{keriso} shows
	that the map $\sum_{u=\ell}^{L}\phi_u : D(\la) \to \D(\mu)$ induces
	a generator of $\Hom_G\bigl(\Delta(\la),\Delta(\mu)\bigr)$.
	
	(2) Suppose $\ell<L$, $\la_1$ is even,
	and $\la_2, \mu_2$ have opposite parity.
	By Lemma \ref{lemmod2}, every vector in
	$\ker_\K\left[\begin{smallmatrix}A\\B\end{smallmatrix}\right]$
	has all coordinates zero except possibly the one indexed by $\ell$.
	Hence $\Hom_G(\Delta(\la),\Delta(\mu)) \neq 0$ if and only if
	column $\ell$ of both matrices is zero. We consider  two cases according to $\ell$.
	
(2a) Let $\ell=0$. Then $\mu_2\le\la_3$ by
\eqref{assumptions2}, and $\mu_2\ge1$ since $\ell<L$.
Thus $(1,1)\in R_B$ by \eqref{RB}, and
\[
B_{(1,1),0}
=\tbinom10\tbinom{\la_2}{0}=1.
\]
Hence $\Hom_G(\Delta(\la),\Delta(\mu))=0$ 
in agreement with the criterion in \eqref{thm(2)}, since $\lambda_3\ge 2^{v_2(\ell+1)}$.

(2b) Let $\ell>0$. Then $\ell=\mu_2-\la_3$.

Let $(t,k)\in R_A$. If $k>\ell$, then
$
\tbinom{\mu_2-k}{\ell-k}=0,
$
and $A_{(t,k),\ell}=0$ by Definition $\ref{defAB}$.
We conclude that column $\ell$ of $A$ is zero if and only if
\[
A_{(t,\ell),\ell}
=\tbinom{\la_1-\ell+t}{t} =0, \ \ 1\le t\le\la_2-\ell.
\]

Similarly, if $(t,k)\in R_B$, then we have
$k\ge\ell+t$ by \eqref{RB}. If $k>\ell+t$, then
$
\tbinom{\la_2+t-k}{t+\ell-k}=0,
$
so $B_{(t,k),\ell}=0$ by Definition $\ref{defAB}$.
We conclude that column $\ell$ of $B$ is zero if and only if
\[
B_{(t,\ell+t),\ell}
=\tbinom{\ell+t}{t} =0, \ \ 1\le t\le\la_3.
\]
By Lemma \ref{lembinom2}, the above conditions
are equivalent to
$
\la_2-\ell<2^{v_2(\la_1-\ell+1)},
\ 
\la_3<2^{v_2(\ell+1)},
$
which is \eqref{thm(2)}. When these inequalities hold, Proposition
\ref{keriso} shows that
$\phi_\ell:D(\la)\to\D(\mu)$ induces a
generator of
$\Hom_G(\Delta(\la),\Delta(\mu))$.

	(3) From Proposition \ref{propiso} and Proposition \ref{keriso} we have the claimed isomorphism and generator of $\Hom_G(\D(\la),\D(\mu))$.
	
	It remains to prove that only  a finite number of halvings can be applied starting with the pair $(\la, \mu)$ . By
	Lemma \ref{lemsmaller}(1),(2), each application of (3)
	produces a pair satisfying \eqref{assumptions1} and
	$\ell'\le L'$, with $0\le\la'_3<\la_3$.
\end{proof}
	
\begin{corollary}\label{dimhom}We have $\dim_\K\Hom_G\bigl(\Delta(\la),\Delta(\mu)\bigr)\le1$ for all dominant weights $\la, \mu$ of $G$.
\end{corollary}
\begin{proof}This follows from the reductions in Section \ref{secred} and Theorem \ref{thm2}.
\end{proof}

\begin{example}
	Let $q\ge2$ be an integer and consider the partitions
	\[
	\la=(8q-6,4q+1,4q-5),
	\ \ 
	\mu=(12q-7,4q-3).
	\]
	For this pair we have
$
	\ell=2, 
	L=4q-3.
$
	Thus $\ell<L$. Since $\la_1$ is even and $\la_2,\mu_2$ are both odd,
	we apply Definition \ref{defhalf} with $\varepsilon=1$. This gives
	the pair $(\la',\mu')$, where
	\[
	\la'=(4q-4,2q,2q-3),
	\ \ 
	\mu'=(6q-5,2q-2).
	\]
	
	For the pair $(\la',\mu')$ we have
$
	\ell'=1, L'=2q-2.
$
	Again $\ell'<L'$. Now $\la'_1$ is even and $\la'_2,\mu'_2$ are both
	even. We apply Definition \ref{defhalf} a second time, but with
	$\varepsilon=0$ to obtain the pair $(\la'',\mu'')$, where
	\[
	\la''=(2q-3,q-1,q-2),
	\ \ 
	\mu''=(3q-4,q-2).
	\]
	
	For $(\la'',\mu'')$ we have $\ell''=0$ and $L''=q-2$.
	Since $\la''_1=2q-3$ is odd, we apply Theorem \ref{thm2}(1)
	to conclude $\Hom_G(\Delta(\la''),\Delta(\mu'')) \neq 0$ if and only if
	$	q-1<2^{v_2(q)}$ and
	$
	q-2<2^{v_2(q)}.
	$
	These inequalities hold if and only if $q$ is a power of 2. By
	Theorem \ref{thm2}(3), we conclude that
	\[
	\dim_\K\Hom_G(\Delta(\la),\Delta(\mu))
	=
	\begin{cases}
		1,&\text{if $q$ is a power of $2$},\\
		0,&\text{otherwise}.
	\end{cases}
	\]
	
	Suppose now that $q$ is a power of $2$. We determine a generator of $\Hom_G(\Delta(\la),\Delta(\mu))$ as follows. Since $\la''=(2q-3,q-1,q-2)$ and $2q-3$ is odd, by
	Theorem \ref{thm2}(1) a generator of $\Hom_G(\D(\la''),\D(\mu''))$
	has coefficients $c''_j=1,  \ (0\le j\le q-2)$.
	
	For the pair $(\la',\mu')$ we have $\varepsilon=0$. By \eqref{inv} a generator of  $\Hom_G(\D(\la'),\D(\mu'))$ has coefficients
$
	c'_{2j+1}=c''_j=1,
	\ 
	c'_{2j+2}=0
	\ \ (0\le j\le q-2).
$
	
	For the pair $(\la,\mu)$ we have $\varepsilon=1$. By \eqref{inv} a generator of  $\Hom_G(\D(\la),\D(\mu))$ has coefficients
	$
	c_{2i}=c_{2i+1}=c'_i
	\  (1\le i\le 2q-2).
	$
	In other words, a generator of
	$\Hom_G(\Delta(\la),\Delta(\mu))$ is induced by the map
$\sum_{j=0}^{q-2}\bigl(\phi_{4j+2}+\phi_{4j+3}\bigr): D(\la) \to \D(\mu)$.
\end{example}

	\subsection{Homomorphisms between Specht modules}\label{secspecht} For a partition $\nu$ of $r$, let $S^\nu$ be the corresponding Specht module of the symmetric group $\mathfrak{S}_r$. Recall that $\nu$ is called $2$-regular if its positive parts are distinct.

	\begin{corollary}\label{specht} Let
		$\la,\mu$ be partitions of an integer $r\ge3$, each with at
		most three parts. Suppose that $\mu$ is $2$-regular, or that
		$\mu_1=\mu_2\ge\mu_3+2$. Then
		\[
		\dim_\K\Hom_{\K\mathfrak S_r}(S^\mu,S^\la)
		=\dim_\K\Hom_{\GL_3(\K)}\bigl(\Delta(\la),\Delta(\mu)\bigr).
		\] 
	\end{corollary}
	We note that $\dim_\K\Hom_{\GL_3(\K)}\bigl(\Delta(\la),\Delta(\mu)\bigr)$ is at most one 
	and is determined by Theorem \ref{thm2}, after the reductions
	of Section \ref{secred}.
	\begin{proof}
		Let $d \ge r$  and let $f$ be the Schur functor from the category of homogeneous polynomial representations of $\GL_d(\K)$ of degree $r$ to the category of  $\K\mathfrak{S}_r$-modules \cite[Section 6.1]{Gr}. Let
		 $G_\otimes$ be the left adjoint of $f$ \cite[Definition 3.3]{McD}. Let $\Delta_d(\nu)$ and
		$\nabla_d(\nu)$ be the Weyl and dual Weyl modules for
		$\GL_d(\K)$. In the second case, the partition $(\mu_2,\mu_3)$ is $2$-regular,
		since $\mu_2\ge\mu_3+2$. Thus, in both cases,
		$G_\otimes(S^\mu)\cong\nabla_d(\mu)$
		by \cite[Theorem 1.2]{McD}. We also have
		$f(\nabla_d(\la))\cong S^\la$ \cite[Section 6.3]{Gr}. Hence 
	$
			\Hom_{\K\mathfrak S_r}(S^\mu,S^\la)
			\cong\Hom_{\GL_d(\K)}
			\bigl(\nabla_d(\mu),\nabla_d(\la)\bigr)
			\cong\Hom_{\GL_d(\K)}
			\bigl(\Delta_d(\la),\Delta_d(\mu)\bigr),
		$
	where we used contravariant duality in the second isomorphism. But for $d \ge3$ and partitions with at most three parts, it is well known that $\Hom_{\GL_d(\K)}
	\bigl(\Delta_d(\la),\Delta_d(\mu)\bigr)$ $ \cong \Hom_{\GL_3(\K)}
	\bigl(\Delta_3(\la),\Delta_3(\mu)\bigr)$, see, for example, the proof of \cite[Proposition 1.1]{Ku}. \end{proof}
\section{Proofs of Lemmas \ref{lemsmaller} and \ref{lemrowhalf}}\label{lemproofs} In this section, we prove Lemma \ref{lemsmaller} and Lemma \ref{lemrowhalf}

We adopt throughout the assumptions and notation of Definition \ref{defhalf}.
\subsection{Proof of Lemma \ref{lemsmaller}}\label{proof33}
(1) By \eqref{assumptions2} and the assumption $\ell<L$,
we have
\[
0\le\ell<L\le\mu_2,\qquad
\mu_2-\la_3\le\ell<L\le\mu_2.
\]
Hence $\mu_2,\la_3\ge1$, so Definition \ref{defhalf}
gives $\mu'_2,\la'_3\ge0$.

We first show that $\lambda'$ is a partition. Since $\la_1$ is even and
$\la_2,\mu_2\le\la_1$, we have
$
\la'_2,\mu'_2
\le\left\lfloor{(\la_1-1)}/{2}\right\rfloor
=\frac{\la_1}{2}-1=\la'_1.
$
Moreover, since $\la_3\le\la_2$ and
$\la_2\equiv\varepsilon\pmod2$, we have
\[
\la'_3
=\left\lfloor{(\la_3-1+\varepsilon)}/{2}\right\rfloor
\le\left\lfloor{(\la_2-1+\varepsilon)}/{2}\right\rfloor
=\left\lfloor{(\la_2-1)}/{2}\right\rfloor
=\la'_2.
\]
Thus $\la'_1\ge\la'_2\ge\la'_3\ge0$,
so $\la'$ is a partition. We also have
$0\le\mu'_2\le\la'_1$.

We now show that $\mu'$ is a partition. By \eqref{assumptions2} and $\ell<L$ we get
$
\mu_2-\la_3\le\ell<L\le\la_2,
$
so $\mu_2-\la_2<\la_3$.
Since $\la_2\equiv\mu_2\pmod2$, the formulas
in Definition \ref{defhalf} give
$
2(\mu'_2-\la'_2)=\mu_2-\la_2.
$
Also, the formula for $\la'_3$ gives
$\la_3\le2\la'_3+2$. Hence
$
2(\mu'_2-\la'_2)
=\mu_2-\la_2
<\la_3\le2\la'_3+2.
$
Since $\mu'_2-\la'_2$ and $\la'_3$ are integers,
we obtain
$\mu'_2-\la'_2\le\la'_3$.
Therefore,
\[
\mu'_1
=\la'_1+\la'_2+\la'_3-\mu'_2
\ge\la'_1\ge\mu'_2
\]
and $\mu'$ is a partition. We have
$|\la'|=|\mu'|$ by the definition of $\mu'_1$.

Finally, since $\la_3\ge1$ and
$\varepsilon\in\{0,1\}$, Definition \ref{defhalf}
gives
$
0\le\la'_3
=\left\lfloor{(\la_3-1+\varepsilon)}/{2}\right\rfloor
\le\frac{\la_3}{2}<\la_3.
$

(2) Since $\mu_2\equiv\varepsilon\pmod 2$, the formulas
for $\mu'_2$ and $\lambda'_3$ in Definition \ref{defhalf} give
$
\mu'_2-\lambda'_3
=\left\lfloor{(\mu_2-\lambda_3)}/{2}\right\rfloor.
$
Thus $\ell'=\max(0,\mu'_2-\lambda'_3)=\lfloor\ell/2\rfloor$.
Similarly, $L'=\lfloor(L-1)/2\rfloor$.
Since $\ell<L$, we have $\ell'\le L'$.

(3) For every integer $i$, part (2) gives
\[
\ell\le2i+1\le L
\iff
\left\lfloor{\ell}/{2}\right\rfloor
\le i\le
\left\lfloor{(L-1)}/{2}\right\rfloor
\iff
\ell'\le i\le L'.
\]

\subsection{Proof of Lemma \ref{lemrowhalf}}\label{proof34}
(1) Let $h\ge1$ and $j\ge0$.

(1A) We first consider $R_A$. Suppose $(2h,2j+1)\in R_A$. By \eqref{RA},
$
\ell\le2j+1\le\min(\la_2-2h,\mu_2).
$
By Lemma \ref{lemsmaller}(2) and
Definition \ref{defhalf}, these inequalities give
\begin{equation}\label{lem341}
\ell'\le j\le\min(\la'_2-h,\mu'_2).
\end{equation}
Thus $(h,j)\in R_{A'}$ by \eqref{RA}.

Conversely, suppose $(h,j)\in R_{A'}$.
Then we have $1\le h\le\la'_2$ and \eqref{lem341}.
Lemma \ref{lemsmaller}(2) and
Definition \ref{defhalf} yield
\[
\ell\le2j+1\le\min(\la_2-2h,\mu_2).
\]
Moreover, Definition \ref{defhalf} gives
$1\le2h\le\la_2$.
Hence $(2h,2j+1)\in R_A$.

(1B) We next consider $R_B$. Suppose $(2h,2j+1)\in R_B$.
By \eqref{RB},
\[
\max(0,\mu_2-\la_3+2h)
\le2j+1\le\min(\la_2+2h,\mu_2).
\]
Using Definition \(\ref{defhalf}\) we see that these inequalities give
\[
\max(0,\mu'_2-\la'_3+h)
\le j\le\min(\la'_2+h,\mu'_2).
\]
In particular,
$\mu'_2-\la'_3+h\le j\le\mu'_2$,
so $1\le h\le\la'_3$.
Thus $(h,j)\in R_{B'}$ by \eqref{RB}.

Conversely, suppose $(h,j)\in R_{B'}$.
Then
\begin{equation}\label{lem342}
	\max(0,\mu'_2-\la'_3+h)
	\le j\le\min(\la'_2+h,\mu'_2).
\end{equation}
By \eqref{lem342} and Definition \ref{defhalf} we get
$
\max(0,\mu_2-\la_3+2h)
\le2j+1\le\min(\la_2+2h,\mu_2).
$
Moreover, Definition \ref{defhalf} gives
$1\le2h\le\la_3$.
Hence $(2h,2j+1)\in R_B$.

(2) Suppose $\varepsilon=1$, so that
$\la_2$ and $\mu_2$ are odd.

Suppose $(2h,2j)\in R_A$. \eqref{RA} implies
$
\ell\le2j\le\min(\la_2-2h,\mu_2).
$
Since this upper bound is odd,
\[
\ell\le2j+1\le\min(\la_2-2h,\mu_2).
\]
Thus
$(2h,2j+1)\in R_A$. By part (1) of the lemma we have
$(h,j)\in R_{A'}$.

Finally, suppose $(2h,2j)\in R_B$.
By \eqref{RB} we have
$
\max(0,\mu_2-\la_3+2h)
\le2j\le\min(\la_2+2h,\mu_2).
$
Agaim the upper bound is odd, and thus 
$
\max(0,\mu_2-\la_3+2h)
\le2j+1\le\min(\la_2+2h,\mu_2).
$
Hence 
$(2h,2j+1)\in R_B$. From part  (1) of the lemma we conclude that
$(h,j)\in R_{B'}$.

\end{document}